\documentclass[leqno]{article}
\usepackage{a4,latexsym,amssymb,exscale}

\begin{document}

\newtheorem{defn}{Definition}[section]
\newtheorem{thm}[defn]{Theorem}
\newtheorem{lemma}[defn]{Lemma}
\newtheorem{prop}[defn]{Proposition}
\newtheorem{cor}[defn]{Corollary}
\newenvironment{proof}{\itshape Proof: \upshape}{\hfill\ensuremath{\Box}\medskip}

\renewcommand{\arraystretch}{1.5}
\renewcommand{\labelenumi}{\roman{enumi})}
\renewcommand{\O}{\mathcal{O}}
\newcommand{\A}{\mathbb{A}}
\renewcommand{\L}{\mathcal{L}}
\newcommand{\E}{\mathcal{E}}
\newcommand{\F}{\mathcal{F}}
\newcommand{\V}{\mathcal{V}}
\newcommand{\Gl}{\mathrm{Gl}}
\newcommand{\Spec}{\mathrm{Spec}}
\newcommand{\discr}{\mathfrak{d}}
\newcommand{\card}{\mathrm{card}}
\renewcommand{\div}{\mathrm{div}}
\newcommand{\dual}{\mathrm{dual}}
\newcommand{\Sl}{\mathrm{Sl}}
\newcommand{\Mat}[3][{}]{\mathrm{Mat}^{#1}_{#2 \times #3}}
\newcommand{\id}{\mathrm{id}}
\newcommand{\rk}{\mathrm{rk}}
\newcommand{\reals}{\mathbb{R}}
\newcommand{\nonnegs}{\reals_{\geq 0}}
\newcommand{\integers}{\mathbb{Z}}
\newcommand{\naturals}{\mathbb{N}}
\newcommand{\rationals}{\mathbb{Q}}
\newcommand{\complexnums}{\mathbb{C}}
\newcommand{\sections}[1]{\Gamma(#1)}
\newcommand{\landau}{\mathrm{o}}
\newcommand{\Landau}{\mathrm{O}}
\newcommand{\eps}{\varepsilon}

\newcommand{\longto}[1][]{\stackrel{#1}{\longrightarrow}}
\newcommand{\blank}{\hspace{.5ex}.\hspace{.5ex}}

\title{Stability of Arakelov bundles and tensor products without global sections}
\author{Norbert Hoffmann}

\maketitle

\section*{Introduction}

G. Faltings has proved that for each semistable vector bundle $E$ over an algebraic curve of genus $g$, there is
another vector bundle $F$ such that $E \otimes F$ has slope $g-1$ and no global sections. (Note that any vector bundle
of slope $> g-1$ has global sections by Riemann--Roch.) See \cite{g_bundles} and \cite{stacks} where this result is
interpreted in terms of theta functions and used for a new construction of moduli schemes of vector bundles.

In the present paper, an arithmetic analogue of that theorem is proposed. The algebraic curve is replaced by the set
$X$ of all places of a number field $K$; we call $X$ an arithmetic curve. Vector bundles are replaced by so-called
Arakelov bundles, cf. section \ref{bundles}. In the special case $K = \rationals$, Arakelov bundles without global
sections are lattice sphere packings, and the slope $\mu$ measures the packing density.

We will see at the end of section \ref{no_sections} that the maximal slope of Arakelov bundles of rank $n$ without
global sections is $d(\log n + \Landau(1) )/2 + (\log \discr)/2$ where $d$ is the degree and $\discr$ is the
discriminant of $K$. Now the main result is:
\begin{thm} \label{thm_intro}
  Let $\E$ be a semistable Arakelov bundle over the arithmetic curve $X$. For each $n \gg 0$ there is an Arakelov
  bundle $\F$ of rank $n$ satisfying
  \begin{displaymath}
    \mu(\E \otimes \F) > \frac{d}{2}(\log n - \log \pi - 1 - \log 2) + \frac{\log \discr}{2}
  \end{displaymath}
  such that $\E \otimes \F$ has no nonzero global sections.
\end{thm}

The proof is inspired by (and generalizes) the Minkowski-Hlawka existence theorem for sphere packings; in particular,
it is not constructive. The principal ingredients are integration over a space of Arakelov bundles (with respect to
some Tamagawa measure) and an adelic version of Siegel's mean value formula. Section \ref{mean_value} explains the
latter, section \ref{bundles} contains all we need about Arakelov bundles, and the main results are proved and
discussed in section \ref{no_sections}.

This paper is a condensed and slightly improved part of the author's Ph.\,D. thesis \cite{diss}. I would like to thank
my adviser G. Faltings for his suggestions; the work is based on his ideas. It was supported by a grant of the
Max-Planck-Institut in Bonn.

\section{Notation} \label{notation}

Let $K$ be a number field of degree $d$ over $\rationals$ and with ring of integers $\O_K$. Let $X = \Spec(\O_K)
\cup X_{\infty}$ be the set of places of $K$; this might be called an `arithmetic curve' in the sense of Arakelov
geometry. $X_{\infty}$ consists of $r_1$ real and $r_2$ complex places with $r_1+2r_2 = d$. $w(K)$ is the number
of roots of unity in $K$.

For every place $v \in X$, we endow the corresponding completion $K_v$ of $K$ with the map $| \blank |_v : K_v \to
\nonnegs$ defined by $\mu(a \cdot S) = |a|_v \cdot \mu(S)$ for a Haar measure $\mu$ on $K_v$. This is the normalized
valuation if $v$ is finite, the usual absolute value if $v$ is real and its square if $v$ is complex. The well known
product formula $\prod_{v \in X} |a|_v = 1$ holds for every $0 \neq a \in K$. On the adele ring $\A$, we have the
divisor map $\div: \A \to \nonnegs^X$ that maps each adele $a = (a_v)_{v \in X}$ to the collection
$(|a_v|_v)_{v \in X}$ of its norms.

Let $\O_v$ be the set of those $a \in K_v$ which satisfy $|a|_v \leq 1$; this is the ring of integers in $K_v$ for
finite $v$ and the unit disc for infinite $v$. Let $\O_{\A}$ denote the product $\prod_{v \in X} \O_v$; this is
the set of all adeles $a$ with $\div(a) \leq 1$. By $D \leq 1$ for an element $D = (D_v)_{v \in X}$ of $\nonnegs^X$,
we always mean $D_v \leq 1$ for all $v$.

We fix a canonical Haar measure $\lambda_v$ on $K_v$ as follows:
\begin{itemize}
 \item If $v$ is finite, we normalize by $\lambda_v(\O_v) = 1$.
 \item If $v$ is real, we take for $\lambda_v$ the usual Lebesgue measure on $\reals$.
 \item If $v$ is complex, we let $\lambda_v$ come from the real volume form $dz \wedge d\bar{z}$ on $\complexnums$.
  In other words, we take twice the usual Lebesgue measure.
\end{itemize}
This gives us a canonical Haar measure $\lambda := \prod_{v \in X} \lambda_v$ on $\A$. We have $\lambda(\A/K) =
\sqrt{\discr}$ where $\discr= \discr_{K/\rationals}$ denotes (the absolute value of) the discriminant. More details
on this measure can be found in \cite{weil}, section 2.1.

Let $V_n = \frac{\pi^{n/2}}{(n/2)!}$ be the volume of the unit ball in $\reals^n$. For $v \in X_{\infty}$, we denote
by $\O_v^n$ the unit ball with respect to the standard scalar product on $K_v^n$. Observe that this is \emph{not} the
$n$-fold cartesian product of $\O_v \subseteq K_v$. Similarly, $\O_{\A}^n := \prod_{v \in X} \O_v^n$ is not the
$n$-fold product of $\O_{\A} \subseteq \A$. Its volume $\lambda^n(\O_{\A}^n)$ is $V_n^{r_1}(2^n V_{2n})^{r_2}$.

\section{A mean value formula} \label{mean_value}

The following proposition is a generalization of Siegel's mean value formula to an adelic setting: With real
numbers and integers instead of adeles and elements of $K$, Siegel has already stated it in \cite{siegel}, and
an elementary proof is given in \cite{macbeath}. (In the special case $l=1$, a similar question is studied in
\cite{thunder}.)

\begin{prop}
  Let $1 \leq l < n$, and let $f$ be a nonnegative measurable function on the space $\Mat{n}{l}(\A)$ of
  $n \times l$ adele matrices. Then
  \begin{equation} \label{mittelwertformel}
    \int\limits_{\Sl_n(\A)/\Sl_n(K)}
      \sum_{\textstyle {M \in \Mat{n}{l}(K) \atop \rk(M)=l}} f(g \cdot M) \,d\tau(g)
    = \discr^{-nl/2} \int\limits_{\Mat{n}{l}(\A)} f\,d\lambda^{n \times l}
  \end{equation}
  where $\tau$ is the unique $\Sl_n(\A)$-invariant probability measure on $\Sl_n(\A)/\Sl_n(K)$.
\end{prop}
\begin{proof}
  The case $l=1$ is done in section 3.4 of \cite{weil}, and the general case can be deduced along the same
  lines from earlier sections of this book. We sketch the main arguments here; more details are given in
  \cite{diss}, section 3.2.

  Let $G$ be the algebraic group $\Sl_n$ over the ground field $K$, and denote by $\tau_G$ the Tamagawa measure
  on $G(\A)$ or any quotient by a discrete subgroup. The two measures $\tau$ and $\tau_G$ on $\Sl_n(\A)/\Sl_n(K)$
  coincide because the Tamagawa number of $G$ is one.

  $G$ acts on the affine space $\Mat{n}{l}$ by left multiplication. Denote the first $l$ columns of the $n \times
  n$ identity matrix by $E \in \Mat{n}{l}(K)$, and let $H \subseteq G$ be the stabilizer of $E$. This algebraic
  group $H$ is a semi-direct product of $\Sl_{n-l}$ and $\Mat{l}{(n-l)}$. Hence section 2.4 of \cite{weil} gives
  us a Tamagawa measure $\tau_H$ on $H(\A)$, and the Tamagawa number of $H$ is also one.

  Again by section 2.4 of \cite{weil}, we have a Tamagawa measure $\tau_{G/H}$ on $G(\A)/H(\A)$ as well, and it
  satisfies $\tau_G = \tau_{G/H} \cdot \tau_H$ in the sense defined there. In particular, this implies
  \begin{displaymath}
    \int\limits_{G(\A)/H(K)} f(g \cdot E)\,d\tau_G(g) = \int\limits_{G(\A)/H(\A)} f(g \cdot E)\,d\tau_{G/H}(g).
  \end{displaymath}
  It is easy to see that the left hand sides of this equation and of (\ref{mittelwertformel}) coincide.
  According to lemma 3.4.1 of \cite{weil}, the right hand sides coincide, too.
\end{proof}

\section{Arakelov vector bundles} \label{bundles}

Recall that a (euclidean) lattice is a free $\integers$-module $\Lambda$ of finite rank together with a scalar
product on $\Lambda \otimes \reals$. This is the special case $K = \rationals$ of the following notion:
\begin{defn} \upshape
  An \emph{Arakelov (vector) bundle} $\E$ over our arithmetic curve $X = \Spec(\O_K) \cup X_{\infty}$ is a finitely
  generated projective $\O_K$-module $\E_{\O_K}$ endowed with
  \begin{itemize}
   \item a euclidean scalar product $\langle \_\,, \_\, \rangle_{\E, v}$ on the real vector space $\E_{K_v}$ for
    every real place $v \in X_{\infty}$ and
   \item an hermitian scalar product $\langle \_\,, \_\, \rangle_{\E, v}$ on the complex vector space $\E_{K_v}$
    for every complex place $v \in X_{\infty}$
  \end{itemize}
  where $\E_A := \E_{O_K} \otimes A$ for every $\O_K$-algebra $A$.
\end{defn}

A first example is the trivial Arakelov line bundle $\O$. More generally, the trivial Arakelov vector bundle $\O^n$
consists of the free module $\O_K^n$ together with the standard scalar products at the infinite places.

We say that $\E'$ is a subbundle of $\E$ and write $\E' \subseteq \E$ if $\E'_{\O_K}$ is a direct summand in
$\E_{\O_K}$ and the scalar product on $\E'_{K_v}$ is the restriction of the one on $\E_{K_v}$ for every infinite
place $v$. Hence every vector subspace of $\E_K$ is the generic fibre of one and only one subbundle of $\E$.

From the data belonging to an Arakelov bundle $\E$, we can define a map
\begin{displaymath}
  \| \blank \|_{\E,v} : \E_{K_v} \longto \nonnegs
\end{displaymath}
for every place $v \in X$:
\begin{itemize}
 \item If $v$ is finite, let $\|e\|_{\E,v}$ be the minimum of the valuations $|a|_v$ of those elements $a \in K_v$
  for which $e$ lies in the subset $a \cdot \E_{\O_v}$ of $\E_{K_v}$. This is the nonarchimedean norm corresponding
  to $\E_{\O_v}$.
 \item If $v$ is real, we put $\|e\|_{\E,v} := \sqrt{\langle e,e \rangle_v}$, so we just take the norm coming from
  the given scalar product.
 \item If $v$ is complex, we put $\|e\|_{\E,v} := \langle e,e \rangle_v$ which is the square of the norm coming from
  our hermitian scalar product.
\end{itemize}
Taken together, they yield a divisor map
\begin{displaymath}
  \div_{\E}: \E_{\A} \to \nonnegs^X \qquad \qquad e = (e_v) \mapsto (\|e_v\|_{\E, v}).
\end{displaymath}
Although $\O_{\A}$ is not an $\O_K$-algebra, we will use the notation $\E_{\O_{\A}}$, namely for the compact set
defined by
\begin{displaymath}
  \E_{\O_{\A}} := \{e \in \E_{\A}: \div_{\E}(e) \leq 1\}.
\end{displaymath}

Recall that these norms are used in the definition of the Arakelov degree: If $\L$ is an Arakelov line bundle and
$0 \neq l \in \L_K$ a nonzero generic section, then
\begin{displaymath}
  \deg(\L) := -\log \prod_{v \in X} \|l\|_{\L,v}
\end{displaymath}
and the degree of an Arakelov vector bundle $\E$ is by definition the degree of the Arakelov line bundle $\det(\E)$.
$\mu(\E) := \deg(\E)/\rk(\E)$ is called the slope of $\E$. One can form the tensor product of two Arakelov bundles
in a natural manner, and it has the property $\mu(\E \otimes \F) = \mu(\E) + \mu(\F)$.

Moreover, the notion of stability is based on slopes: For $1 \leq l \leq \rk(\E)$, denote by $\mu^{(l)}_{\max}$
the supremum (in fact it is the maximum) of the slopes $\mu(\E')$ of subbundles $\E' \subseteq \E$ of rank $l$.
$\E$ is said to be stable if $\mu^{(l)}_{\max} < \mu(\E)$ holds for all $l < \rk(\E)$, and semistable if
$\mu^{(l)}_{\max} \leq \mu(\E)$ for all $l$.

To each projective variety over $K$ endowed with a metrized line bundle, one can associate a zeta function as in
\cite{franke} or \cite{manin}. We recall its definition in the special case of Grassmannians associated to Arakelov
bundles:
\begin{defn} \upshape
  If $\E$ is an Arakelov bundle over $X$ and $l \leq \rk(\E)$ is a positive integer, then we define
  \begin{displaymath}
    \zeta_{\E}^{(l)}(s) := \sum_{\textstyle{\E' \subseteq \E \atop \rk(\E') = l}} \exp(s \cdot \deg(\E')).
  \end{displaymath}
\end{defn}
The growth of these zeta functions is related to the stability of $\E$. More precisely, we have the following
asymptotic bound:
\begin{lemma} \label{zetaasymp}
  There is a constant $C = C(\E)$ such that
  \begin{displaymath}
    \zeta_{\E}^{(l)}(s) \leq C \cdot \exp(s \cdot l\mu_{\max}^{(l)}(\E))
  \end{displaymath}
  for all sufficiently large real numbers $s$.
\end{lemma}
\begin{proof}
  Fix $\E$ and $l$. Denote by $N(T)$ the number of subbundles $\E' \subseteq \E$ of rank $l$ and degree at least $-T$.
  There are $C_1, C_2 \in \reals$ such that
  \begin{displaymath}
    N(T) \leq \exp(C_1 T + C_2)
  \end{displaymath}
  holds for all $T \in \reals$. (Embedding the Grassmannian into a projective space, this follows easily from
  \cite{schanuel}. See \cite{diss}, lemma 3.4.8 for more details.)

  If we order the summands of $\zeta_{\E}^{(l)}$ according to their magnitude, we thus get
  \begin{eqnarray*}
    \zeta_{\E}^{(l)}(s) & \leq & \sum_{\nu = 0}^{\infty}
      N \left( -l\mu_{\max}^{(l)}(\E) + \nu + 1 \right) \cdot \exp\left(s \cdot (l\mu_{\max}^{(l)}(\E) - \nu) \right)\\
    & \leq & \exp(s \cdot l\mu_{\max}^{(l)}(\E)) \cdot \sum_{\nu = 0}^{\infty} \frac{C_3}{\exp((s-C_1) \nu)}.
  \end{eqnarray*}
  But the last sum is a convergent geometric series for all $s > C_1$ and decreases as $s$ grows, so it is bounded
  for $s \geq C_1 + 1$.
\end{proof}

\section{The main theorem} \label{no_sections}

The global sections of an Arakelov bundle $\E$ over $X = \Spec(\O_K) \cup X_{\infty}$ are by definition the
elements of the finite set
\begin{displaymath}
  \sections{\E} := \E_K \cap \E_{\O_{\A}} \subseteq \E_{\A}.
\end{displaymath}

Note that in the special case $K = \rationals$, an Arakelov bundle without nonzero global sections is nothing but
a (lattice) sphere packing: $\sections{\E} = 0$ means that the (closed) balls of radius $1/2$ centered at the
points of the lattice $\E_{\integers}$ are disjoint. Here larger degree corresponds to denser packings.

\begin{thm} \label{mainthm}
  Let $\E$ be an Arakelov bundle over the arithmetic curve $X$. If an integer $n > \rk(\E)$ and an Arakelov 
  line bundle $\L$ satisfy
  \begin{displaymath}
    1 > \sum_{l=1}^{\rk(\E)} 
      \frac{1}{\discr^{nl/2}}
      \cdot \lambda^{nl} \left( \frac{K^* \O_{\A}^{nl}}{K^*} \right)
      \cdot \zeta_{\E}^{(l)}(n) \exp(l \deg(\L)),
  \end{displaymath}
  then there is an Arakelov bundle $\F$ of rank $n$ and determinant $\L$ such that
  \begin{displaymath}
    \sections{\E \otimes \F} = 0.
  \end{displaymath}
\end{thm}
\begin{proof}
  Note that any global section of $\E \otimes \F$ is already a global section of $\E' \otimes \F$ for a unique
  minimal subbundle $\E' \subseteq \E$, namely the subbundle whose generic fibre is the image of the induced map
  $(\F_K)^{\dual} \to \E_K$. We are going to average the number of these sections (up to $K^*$) for a fixed
  subbundle $\E'$ of rank $l$.

  Fix one particular Arakelov bundle $\F$ of rank $n$ and determinant $\L$. Choose linear isomorphisms
  $\phi_{\E'}: K^l \to \E_K'$ and $\phi_{\F}: K^n \to \F_K$ and let
  \begin{displaymath}
    \phi: \Mat{n}{l}(K) \longto[\sim] (\E' \otimes \F)_K
  \end{displaymath}
  be their tensor product. Our notation will not distinguish these maps from their canonical extensions to
  completions or adeles.

  For each $g \in \Sl_n(\A)$, we denote by $g \F$ the Arakelov bundle corresponding to the $K$-lattice
  $\phi_{\F}(g K^n) \subseteq \F_{\A}$. More precisely, $g \F$ is the unique Arakelov bundle satisfying
  $(g \F)_{\A} = \F_{\A}$, $(g \F)_{\O_{\A}} = \F_{\O_{\A}}$ and $(g \F)_{K} = \phi_{\F}(g K^n)$. This gives
  the usual identification between $\Sl_n(\A)/\Sl_n(K)$ and the space of Arakelov bundles of rank $n$ and
  fixed determinant together with local trivialisations.

  Observe that the generic fibre of $\E' \otimes g \F$ is $\phi(g \Mat{n}{l}(K))$. A generic section is not in
  $\E'' \otimes g \F$ for any $\E'' \subsetneq \E'$ if and only if the corresponding matrix has rank $l$. So
  according to the mean value formula of section \ref{mean_value}, the average number of global sections
  \begin{displaymath}
    \int\limits_{\Sl_n(\A)/\Sl_n(K)} \card\left(
      \frac{K^* \sections{\E'  \otimes g \F}}{K^*} \setminus \bigcup_{\E'' \subsetneq \E'}
      \frac{K^* \sections{\E'' \otimes g \F}}{K^*}
    \right) \,d\tau(g)
  \end{displaymath}
  is equal to the integral
  \begin{equation} \label{integral}
    \discr^{-nl/2} \int\limits_{\Mat{n}{l}(\A)}
      (f_K \circ \div_{\E' \otimes \F} \circ \phi) \,d\lambda^{n \times l}.
  \end{equation}
  Here the function $f_K: \nonnegs^X \to \nonnegs$ is defined by
  \begin{displaymath}
    f_K(D) := \left\{ \begin{array}{ll}
      1/\card\{a \in K^*: \div(a) \cdot D \leq 1\} & \mbox{if $D \leq 1$}\\
      0                                            & \mbox{otherwise}
    \end{array} \right.
  \end{displaymath}
  with the convention $1/\infty = 0$.

  In order to compute (\ref{integral}), we start with the local transformation formula
  \begin{displaymath} \setlength{\arraycolsep}{2pt} \begin{array}{cccccccc}
    \lambda_v^{n\times l} & (\{M\in & \Mat{n}{l}(K_v) & :c_1\leq & \|\phi(M)\|_{\E'\otimes\F,v} & \leq c_2\}) & = &\\
    \lambda_v^{nl}        & (\{M\in &    K_v^{nl}     & :c_1\leq &    \| M \|                   & \leq c_2\}) &
      \cdot & \| \det(\phi) \|^{-1}_{\det( \E' \otimes \F), v}
  \end{array} \end{displaymath}
  for all $c_1, c_2 \in \nonnegs$. Regarding this as a relation between measures on $\nonnegs$ and taking the
  product over all places $v \in X$, we get the equation
  \begin{equation} \label{trafo}
    (\div_{\E' \otimes \F} \circ \phi)_* \lambda^{n \times l}
    = \exp \deg( \E' \otimes \F ) \cdot (\div_{\O^{nl}})_* \lambda^{nl}
  \end{equation}
  of measures on $\nonnegs^X$. Hence the integrals of $f_K$ with respect to these measures also coincide:
  \begin{displaymath}
    \int\limits_{\Mat{n}{l}(\A)} (f_K \circ \div_{\E' \otimes \F} \circ \phi) \,d\lambda^{n \times l}
    = \exp( n \deg(\E') + l \deg(\F) ) \cdot \lambda^{nl} \left( \frac{K^* \O_{\A}^{nl}}{K^*} \right).
  \end{displaymath}

  We substitute this for the integral in (\ref{integral}). A summation over all nonzero subbundles $\E' \subseteq \E$
  yields
  \begin{displaymath}
    \int\limits_{\Sl_n(\A)/\Sl_n(K)} \hspace{-5.3ex}
      \card\left( \frac{K^* \sections{\E \otimes g \F} \setminus 0}{K^*} \right)
    \,d\tau(g) = \sum_{l=1}^{\rk(\E)} 
      \frac{\zeta_{\E}^{(l)}(n) \exp(l \deg(\F))}{\discr^{nl/2}}
      \lambda^{nl} \left( \frac{K^* \O_{\A}^{nl}}{K^*} \right).
  \end{displaymath}
  But the right hand side was assumed to be less than one, so there there has to be a $g \in \Sl_n(\A)$ with
  $\sections{\E \otimes g \F} = 0$.
\end{proof}

In order to apply this theorem, one needs to compute $\lambda^{N} (K^* \O_{\A}^{N}/K^*)$ for $N \geq 2$. We start
with the special case $K = \rationals$. Here each adele $a \in \O_{\A}^N$ outside a set of measure zero has a
rational multiple in $\O_{\A}^N$ with norm $1$ at all finite places, and this multiple is unique up to sign.
Hence we conclude
\begin{displaymath}
  \lambda^{N} \left( \frac{\rationals^* \O_{\A}^{N}}{\rationals^*} \right)
  = \frac{V_N}{2} \cdot \prod_{p \mbox{ \scriptsize prime}} \lambda_p^N( \integers_p^N \setminus p\integers_p^N)
  = \frac{V_N}{2 \zeta(N)}.
\end{displaymath}
In particular, the special case $K = \rationals$ and $\E = \O$ of the theorem above is precisely the Minkowski-Hlawka
existence theorem for sphere packings \cite{minkowski}, \S15.

For a general number field $K$, we note that the roots of unity preserve $\O_{\A}^N$. Then we apply Stirling's
formula to the factorials occurring via the unit ball volumes and get
\begin{displaymath}
  \lambda^{N} \left( \frac{K^* \O_{\A}^{N}}{K^*} \right) \leq \frac{\lambda^N( \O_{\A}^N)}{w(K)}
  \leq \left( \frac{2 \pi e}{N} \right)^{dN/2} \cdot \left( \frac{1}{\pi N} \right)^{(r_1 + r_2)/2}
    \cdot \frac{1}{2^{r_2/2} w(K)}.
\end{displaymath}
Using such a bound and the asymptotic statement \ref{zetaasymp} about $\zeta^{(l)}_{\E}$, one can deduce the
following corollary of theorem \ref{mainthm}.
\begin{cor} \label{maincor}
  Let the Arakelov bundle $\E$ over $X$ be given. If $n$ is a sufficiently large integer and $\mu$ is a real number
  satisfying
  \begin{displaymath}
    \mu^{(l)}_{\max}(\E) + \mu \leq \frac{d}{2} (\log n + \log l - \log \pi - 1 - \log 2) + \frac{\log \discr}{2}
  \end{displaymath}
  for all $1 \leq l \leq \rk(\E)$, then there is an Arakelov bundle $\F$ of rank $n$ and slope larger than $\mu$ such
  that $\sections{\E \otimes \F} = 0$.
\end{cor}

If $\E$ is semistable, this gives the theorem \ref{thm_intro} stated in the introduction. Here is some evidence that
these bounds are not too far from being optimal:
\begin{prop} \label{converse}
  Assume given $\epsilon > 0$ and a nonzero Arakelov bundle $\E$. Let $n > n(\epsilon)$ be a sufficiently large
  integer, and let $\mu$ be a real number such that
  \begin{displaymath}
    \mu^{(l)}_{\max}(\E) + \mu
      \geq \frac{d}{2} (\log n + \log l - \log \pi - 1 + \log 2 + \epsilon) + \frac{\log \discr}{2}
  \end{displaymath} 
  holds for at least one integer $1 \leq l \leq \rk(\E)$. Then there is no Arakelov bundle $\F$ of rank $n$ and
  slope $\mu$ with $\Gamma( \E \otimes \F) = 0$.
\end{prop}
\begin{proof}
  Fix such an $l$ and a subbundle $\E' \subseteq \E$ of rank $l$ and slope $\mu^{(l)}_{\max}(\E)$. For each $\F$ of
  rank $n$ and slope $\mu$, we consider the Arakelov bundle $\F' := \E' \otimes \F$ of rank $nl$. By Stirling's
  formula, the hypotheses on $n$ and $\mu$ imply
  \begin{displaymath}
    \exp \deg (\F') \cdot \lambda^{nl}(\O^{nl}_{\A}) > 2^{nld} \cdot \discr^{nl/2}.
  \end{displaymath}
  Now choose a $K$-linear isomorphism $\phi: K^{nl} \longto[\sim] \F_K'$ and extend it to adeles. Applying the global
  transformation formula (\ref{trafo}), we get
  \begin{displaymath}
    \lambda^{nl}( \phi^{-1} \F_{\O_{\A}}') > 2^{nld} \cdot \lambda^{nl}( \A^{nl}/K^{nl}).
  \end{displaymath}
  According to Minkowski's theorem on lattice points in convex sets (in an adelic version like \cite{thunder},
  theorem 3), $\phi^{-1} (\F_{\O_{\A}}') \cap K^{nl} \neq \{0\}$ follows. This means that $\F'$ --- and hence
  $\E \otimes \F$ --- must have a nonzero global section.
\end{proof}

Observe that the lower bound \ref{maincor} and the upper bound \ref{converse} differ only by the constant $d \log 2$.
So up to this constant, the maximal slope of such tensor products without global sections is determined by the
stability of $\E$, more precisely by the $\mu^{(l)}_{\max}(\E)$.

Taking $E = \O$, we get lower and upper bounds for the maximal slope of Arakelov bundles without global sections, as
mentioned in the introduction. In the special case $\E = \O$ and $K = \rationals$ of lattice sphere packings,
\cite{conway} states that no essential improvement of corollary \ref{maincor} is known whereas several people have
improved the other bound \ref{converse} by constants.

\end{document}